\documentclass[a4paper]{article}

\usepackage{amssymb}

\newtheorem{thm}{Theorem}
\newtheorem{lemma}[thm]{Lemma}

\newtheorem{corollary}{Corollary}
\newtheorem{example}{Example}
\newtheorem{definition}{Definition}
\newtheorem{algorithm}[thm]{Algorithm}

\begin{document}

\title{On the number of mutually disjoint pairs of S-permutation matrices}

\author{Krasimir Yordzhev}
\date{\empty}

\maketitle
\begin{center}
Faculty of Mathematics and Natural Sciences\\
South-West University, Blagoevgrad, Bulgaria\\
             E-mail: yordzhev@swu.bg
\end{center}

\begin{abstract}
This work examines the concept of S-permutation matrices, namely $n^2 \times n^2$   permutation matrices  containing  a single  1 in each canonical $n \times n$ subsquare (block). The article suggests a formula for counting mutually disjoint pairs of $n^2 \times n^2$ S-permutation matrices in the general case by restricting this task to the problem of finding some numerical characteristics of the elements of specially defined for this purpose factor-set of the set of $n \times n$ binary matrices. The paper describe an algorithm that solves the main problem. To do that, every $n\times n$ binary matrix is represented uniquely as a n-tuple of integers.
\end{abstract}

Keyword: binary matrix; S-permutation matrix; disjoint matrices; Sudoku; factor-set; n-tuple of integers

2010 Mathematics Subject Classification: 05B20

\section{Introduction  and notation}

Let $n$ be a positive integer. By $[n]$ we denote the set $[n] =\left\{ 1,2,\ldots ,n\right\}$.

A \emph{binary} (or   \emph{boolean}, or (0,1)-\emph{matrix})  is a matrix all of whose elements belong to the set $\mathfrak{B} =\{ 0,1 \}$. In this paper we will consider only square binary matrices. With $\mathfrak{B}_n$ we will denote the set of all  $n \times n$  binary matrices. With $\mathfrak{B}_{n,k}$ we will denote the set of all   $n\times n$ binary matrices containing exactly  $k$ elements equal to 1.

Two $n\times n$ binary   matrices $A=(a_{ij} )\in \mathfrak{B}_{n}$ and $B=( b_{ij} )\in \mathfrak{B}_{n}$ will be called \emph{disjoint} if there are not integers $i,j\in [n]$ such that $a_{ij} =b_{ij} =1$, i.e. if $a_{ij} =1$ then $b_{ij} =0$ and if $b_{ij} =1$ then $a_{ij} =0$.

Let $n$ be a positive integer and let $A\in \mathfrak{B}_{n^2}$ be a $n^2 \times n^2$ binary matrix.  With the help of $n - 1$ horizontal lines and $n - 1$ vertical lines $A$ has been separated into $n^2$ of number non-intersecting $n\times n$ square sub-matrices $A_{kl}$, $1\le k,l\le n$, e.i.
\begin{equation}\label{matrA}
A =
\left[
\begin{array}{cccc}
A_{11} & A_{12} & \cdots & A_{1n} \\
A_{21} & A_{22} & \cdots & A_{2n} \\
\vdots & \vdots & \ddots & \vdots \\
A_{n1} & A_{n2} & \cdots & A_{nn}
\end{array}
\right] .
\end{equation}

The sub-matrices $A_{kl}$, $1\le k,l\le n$ will be called \emph{blocks}.

A matrix $A\in \mathfrak{B}_{n^2}$ is called an \emph{S-permutation} if in each row, in each column, and in each block of $A$ there is exactly one 1. Let the set of all $n^2 \times n^2$ S-permutation matrices be denoted by $\Sigma_{n^2}$.

The concept of S-permutation matrix was introduced by Geir Dahl  \cite{dahl} in relation to the popular Sudoku puzzle. Sudoku is a very popular game. On the other hand, it is well known that Sudoku matrices are special cases of Latin squares in the class of gerechte designs \cite{Bailey}.

Obviously a square $n^2 \times n^2$ matrix $M$ with elements of $[n^2 ] =\{ 1,2,\ldots ,n^2 \}$ is a Sudoku matrix if and only if there are  matrices $A_1 ,A_2 ,\ldots ,A_{n^2} \in\Sigma_{n^2}$, each two of them are disjoint and such that $P$ can be given in the following way:
\begin{equation}\label{disj}
M=1\cdot A_1 +2\cdot A_2 +\cdots +n^2 \cdot A_{n^2}
\end{equation}

Some algorithms for obtaining random Sudoku matrices and their valuation are described in detail in  \cite{yordzhev_random} and  \cite{Fontana}.

In \cite{Fontana} Roberto Fontana offers an algorithm which randomly gets a family of
$n^2 \times n^2$ mutually disjoint S-permutation matrices, where $n = 2, 3$. In $n = 3$ he
ran the algorithm 1000 times and found 105 different families of nine mutually
disjoint S-permutation matrices. Then using (\ref{disj}) he obtained $9! \cdot 105 = 38\; 102\; 400$ Sudoku matrices.

But it is known~\cite{Felgenhauer} that the total number of $9\times 9$ Sudoku matrices  is
$$9! \cdot 72^2 \cdot 2^7 \cdot 27\; 704\; 267\; 971 = 6\; 670\; 903\; 752\; 021\; 072\; 936\; 960 $$

Thus, in relation with Fontana's algorithm, it looks useful to calculate the probability of two randomly generated S-permutation matrices to be disjoint.
So the question of enumerating all disjoint pairs of  S-permutation matrices naturally arises. This work is devoted to this task.

As we have shown in  \cite{Yordzhev20151}, with hand calculations of the so assigned task with small values of $n$ ($n=2,3$), it is convenient to use the apparatus of graph theory. Unfortunately, when $n\ge 4$ this approach is inefficient. In this article, we will use only the operations of matrix analysis, which are not difficult to process with computers.

\section{A representation of S-permutation matrices}\label{kikiki1}

 Let $n$ be a positive integer. If $z_1 \; z_2 \; \ldots \; z_n$ is a permutation of the elements of the set $[n] =\left\{ 1,2,\ldots ,n\right\}$ and let us shortly denote $\sigma$ this permutation. Then in this case we will denote by $\sigma (i)$ the $i$-th element of this permutation, i.e. $\sigma (i) =z_i$, $i=1,2,\ldots ,n$.

\begin{definition}\label{defPin}
Let $\Pi_n$ denotes the set of all $n\times n$ matrices, constructed such that $\pi\in\Pi_n$ if and only if the following three conditions are true:

{\bf i)} the elements of $\pi$ are ordered pairs of numbers $\langle i,j\rangle$, where $1\le i,j\le n$;

{\bf ii)} if
$$\left[ \langle a_1 , b_1 \rangle \quad \langle a_2 ,b_2 \rangle \quad \cdots  \quad \langle a_n ,b_n \rangle \right]$$
is the $i$-th row of $\pi$ for any $i\in [n] =\{ 1,2,\ldots ,n\}$, then  $a_1 \; a_2 \; \ldots \; a_n$ in this order is a permutation of the elements of the set $[n]$;

{\bf iii)} if
$$\left[
\begin{array}{c}
\langle a_1 ,b_1 \rangle \\
\langle a_2 ,b_2 \rangle \\
\vdots \\
\langle a_n ,b_n \rangle \\
\end{array}
\right]
$$
is the $j$-th column of $\pi$ for any $j\in [n]$, then  $b_1 ,b_2 ,\ldots , b_n$ in this order is a permutation of the elements of the set $[n]$.
\end{definition}

From Definition \ref{defPin}, it follows that  we can represent each row and each column of a matrix $M\in\Pi_n$ with the help of a permutation of elements of the set $[n]$.

Conversely for every $(2n)$-tuple $$\langle \langle \rho_1 ,\rho_2 ,\ldots ,\rho_n \rangle ,\langle \sigma_1 ,\sigma_2 ,\ldots , \sigma_n \rangle \rangle,$$ where
$$\rho_i = \rho_i (1)\; \rho_i (2) \; \ldots \; \rho_i (n),\quad 1\le i\le n$$
$$\sigma_j = \sigma_j (1)\; \sigma_j (2)\; \ldots \; \sigma_j (n),\quad 1\le j\le n$$
are $2n$ permutations of elements of $[n]$ (not necessarily different), then the matrix
$$
\pi =
\left[
\begin{array}{cccc}
\langle \rho_1 (1),\sigma_1 (1)\rangle & \langle \rho_1 (2),\sigma_2 (1)\rangle & \cdots & \langle \rho_1 (n),\sigma_n (1)\rangle \\
\langle \rho_2 (1),\sigma_1 (2)\rangle & \langle \rho_2 (2),\sigma_2 (2)\rangle & \cdots & \langle \rho_2 (n),\sigma_n (2)\rangle \\
\vdots & \vdots & \ddots & \vdots \\
\langle \rho_n (1),\sigma_1 (n)\rangle  & \langle \rho_n (2),\sigma_2 (n)\rangle & \cdots & \langle \rho_n (n),\sigma_n (n)\rangle
\end{array}
\right]
$$
is matrix of $\Pi_n$. Hence
\begin{equation}\label{|Pin|}
\left| \Pi_n \right| =\left( n! \right)^{2n}
\end{equation}

\begin{definition}
We say that matrices $\pi ' =\left[ {p'}_{ij} \right]_{n\times n} \in\Pi_n$ and $\pi '' =\left[ {p''}_{ij} \right]_{n\times n} \in\Pi_n$ are \emph{disjoint}, if ${p'}_{ij} \ne {p''}_{ij}$ for every $i,j\in[n]$.
\end{definition}

\begin{definition}
Let $\pi ' ,\pi '' \in\Pi_n$, $\pi ' =\left[ {p'}_{ij} \right]_{n\times n}$, $\pi '' =\left[ {p''}_{ij} \right]_{n\times n}$ and let  the integers $i,j\in[n]$ are such that ${p'}_{ij} = {p''}_{ij}$. In this case we will say that   ${p'}_{ij}$ and ${p''}_{ij}$ are \emph{component-wise equal  elements}.
\end{definition}

Obviously two $\Pi_n$-matrices are disjoint if and only if they do not have component-wise equal elements.

\begin{example}\label{ex1}
\rm We consider the following $\Pi_3$-matrices:
\end{example}

$$
\pi' =\left[ p_{ij}' \right] =
\left[
\begin{array}{ccc}
\langle 3,1\rangle & \langle 2,1\rangle & \langle 1,2\rangle \\
\langle 2,3\rangle & \langle 3,2\rangle & \langle 1,1\rangle \\
\langle 3,2\rangle & \langle 1,3\rangle & \langle 2,3\rangle
\end{array}
\right]
$$

$$
\pi'' =\left[ p_{ij}'' \right] =
\left[
\begin{array}{ccc}
\langle 3,2\rangle & \langle 1,3\rangle & \langle 2,1\rangle \\
\langle 3,3\rangle & \langle 1,1\rangle & \langle 2,2\rangle \\
\langle 2,1\rangle & \langle 1,2\rangle & \langle 3,3\rangle
\end{array}
\right]
$$

$$
\pi''' =\left[ p_{ij}''' \right] =
\left[
\begin{array}{ccc}
\langle 3,1\rangle & \langle 1,3\rangle & \langle 2,2\rangle \\
\langle 2,2\rangle & \langle 3,1\rangle & \langle 1,1\rangle \\
\langle 2,3\rangle & \langle 1,2\rangle & \langle 3,3\rangle
\end{array}
\right]
$$

Matrices $\pi'$ and $\pi''$ are disjoint, because they do not have component-wise equal elements.

Matrices $\pi'$ and $\pi'''$ are not disjoint, because they have two component-wise equal elements: $p_{11}' =p_{11}''' =\langle 3,1\rangle$ and $p_{23}' =p_{23}''' =\langle 1,1\rangle$.

Matrices $\pi''$ and $\pi'''$ are not disjoint, because they have three component-wise equal elements: $p_{12}'' =p_{12}''' =\langle 1,3\rangle$, $p_{32}'' =p_{32}''' =\langle 1,2\rangle$, and $p_{33}' =p_{33}''' =\langle 3,3\rangle$.

The relationship between S-permutation matrices and the matrices from the set $\Pi_n$ are given by the following theorem:

\begin{thm}\label{l2fhgg}
Let $n$ be an integer, $n\ge 2$. Then there is one to one correspondence between the sets $\Sigma_{n^2}$ and $\Pi_n$.
\end{thm}

Proof. Let $A\in \Sigma_{n^2}$. Then $A$ is constructed with the help of formula (\ref{matrA}) and for every $i,j\in [n]$ in the block $A_{ij} $ there is only one 1 and let this 1 has coordinates $(a_i ,b_j )$. For every $i,j\in [n]$ we obtain ordered pairs of numbers $\langle a_i ,b_j \rangle$ corresponding to these coordinates. As in every row and every column of $A$ there is only one 1, then the matrix $\left[ \alpha_{ij} \right]_{n\times n}$, where $\alpha_{ij} =\langle a_i ,b_j \rangle $, $1\le i,j\le n$, which is obtained by the ordered pairs of numbers is matrix of $\Pi_n$, i.e. matrix for which the conditions i), ii) and iii) are true.

Conversely, let $\left[ \alpha_{ij} \right]_{n\times n} \in \Pi_n$, where $\alpha_{ij} =\langle a_i ,b_j \rangle $, $i,j \in [n]$, $a_i ,b_j \in [n]$. Then for every $i,j\in [n]$ we construct a binary $n\times n$ matrices $A_{ij}$ with only one 1 with coordinates $(a_i ,b_j )$. Then we obtain the matrix of type (\ref{matrA}). According to the properties i), ii) and iii), it is obvious that the obtained matrix is S-permutation matrix.
\hfill $\Box$

\begin{corollary}
The number of all  pairs of disjoint matrices from $\Sigma_{n^2}$ is equal to the number of all pairs of disjoint matrices from $\Pi_n$.
\end{corollary}

Proof. It is easy to see that with respect of the described in Theorem \ref{l2fhgg} one to one correspondence, every pair of disjoint matrices of $\Sigma_{n^2}$ will correspond to a pair of disjoint matrices of $\Pi_n$ and conversely every pair of disjoint matrices of $\Pi_n$ will correspond to a pair of disjoint matrices of $\Sigma_{n^2}$.
\hfill $\Box$

\begin{corollary}   {\rm \cite{dahl}}
The number of all $n^2 \times n^2 $ S-permutation matrices is equal to
\begin{equation}\label{fcrl2}
\left| \Sigma_{n^2} \right| = \left( n! \right)^{2n}
\end{equation}
\end{corollary}

Proof. It follows immediately from Theorem \ref{l2fhgg} and formula (\ref{|Pin|}).
\hfill $\Box$

\section{A formula for counting all disjoint pairs of $n^2 \times n^2$ S-permutation matrices}\label{kikiki2}

Let $A =[a_{ij} ]_{n\times n}\in \mathfrak{B}_n$. We define the following numerical characteristics of the binary matrix A:
\begin{description}
  \item[$r_k (A)$]  -- the number of rows in $A$ having exactly $k$ units, $k=0,1,2,\ldots ,n$;
  \item[$c_k (A)$]  -- the number of columns in $A$ having exactly $k$ units, $k=0,1,2,\ldots ,n$;
  \item[$\psi_k (A)$] = $\displaystyle r_k (A)+c_k (A)$, $k=0,1,2,\ldots ,n$;
  \item[$\varepsilon (A)$] -- the number of units in $A$.
\end{description}

Let $A,B\in \mathfrak{B}_n$. We will say that $A\sim B $, if  $B$ is obtained from  $A$ after dislocation of some of the rows of  $A$. Obviously,  the relation  defined like that is an equivalence relation.  The factor-set ${\mathfrak{B}_n}_{/_\sim}$, i.e. the set of equivalence classes on the above defined relation we denote with $\overline{\mathfrak{B}}_n$. If $A\in \mathfrak{B}_n$, then with $\overline{A}$ we will denote the set $\overline{A} =\{ B\in \mathfrak{B}_n \; |\; B\sim A\}$. Thus $|\overline{A}| =|\{ B\in \mathfrak{B}_n \; |\; B\sim A\} |$ is the cardinality of the set $\overline{A}$. By definition $\overline{\mathfrak{B}}_{n,k} ={\mathfrak{B}_{n,k}}_{/_\sim}$

Obviously if $A,B\in \mathfrak{B}_n$ and $A\sim B$, then $r_k (A)=r_k (B)$, $c_k (A)=c_k (B)$, $\psi_k (A)=\psi_k (B)$, $\varepsilon_k (A)=\varepsilon_k (B)$, $k=0,1,2,\ldots ,n$.  So in a natural way we can define the functions $r_k$, $c_k$, $\psi_k$ and $\varepsilon$ in the factor-set $\overline{\mathfrak{B}}_n ={\mathfrak{B}_n}_{/_\sim}$ as $r_k (\overline{A})$, $c_k (\overline{A})$, $\psi_k (\overline{A})$ and $\varepsilon (\overline{A})$ will mean respectively $r_k (A)$, $c_k (A)$, $\psi_k (A)$ $\varepsilon (A)$, where $A$ is an arbitrary representative of the set $\overline{A} =\{ B\in \mathfrak{B}_n \; |\; B\sim A\}$.

\begin{lemma}\label{l3dskmat}
Let  $\pi\in\Pi_n$. Then the number $q(n,k)$ of all matrices $\pi'  \in\Pi_n$ (including $\pi$),  having at least  $k$, $k=0,1,\ldots ,n^2$ component-wise equal elements to the matrix  $\pi$ is equal to
\begin{equation}\label{fl3gfgfg}
q(n,k)= \sum_{\overline{A}\in \overline{\mathfrak{B}}_{n,k} } |\overline{A} |  \prod_{i=0}^{n-2} \left[ \left( n-i\right) ! \right]^{\psi_i (\overline{A})}
\end{equation}
\end{lemma}

Proof. Let  $\pi =\left[ p_{ij} \right]_{n\times n} ,\pi' =\left[ p'_{ij} \right]_{n\times n} \in \Pi_n$  and let $\pi$  and $\pi'$  have exactly $k$ component-wise equal elements. Then we uniquely  obtain the binary $n\times n$  matrix  $A=\left[ a_{ij} \right]_{n\times n}$, such that $a_{ij} =1$  if and only if $p_{ij} =p'_{ij}$, $i,j\in [n]$.

Inversely, let $A=[a_{ij} ]_{n\times n} \in \mathfrak{B}_n$  and let $\pi =\left[ p_{ij} \right]_{n\times n}$  be an arbitrary matrix from  $\Pi_n$.We search for the number $h(\pi ,A)$ of all matrices  $\pi' =[p'_{ij} ]_{n\times n} \in \Pi_n$, such that  $p'_{ij} =p_{ij}$, if  $a_{ij}=1$. (It is assumed that there exist $s,t\in [n]$  such that $a_{st} =0$  and  $p'_{st} =p_{st}$.)

Let us denote with $\gamma_s$ the number of 1 in $s$-th row of $A$ and let the  $s$-th row of $\pi$  correspond to the permutation  $\rho_s$ of the elements of $[n]$, $s=1,2,\ldots ,n$.  Then there exist $(n-\gamma_s)!$  permutations $\rho'$  of the elements of  $ [n]$, such that if   $a_{st} =1$, then $\rho_s (t) =\rho' (t)$, $t\in [n]$. Likewise we also prove the respective statement for the columns of   $\pi$. Therefore
$$\displaystyle h(\pi ,A) = \prod_{i=0}^n \left[ (n-i)! \right]^{r_i (A)} \prod_{i=0}^n \left[ (n-i)! \right]^{c_i (A)} = \prod_{i=0}^n \left[ (n-i)! \right]^{\psi_i (A)} .$$

From everything said so far it follows that for each $\pi\in\Pi_n$  there exist
$$q(n,k)=\sum_{A\in\mathfrak{B}_{n,k}}  \prod_{i=0}^n \left[ (n-i)! \right]^{\psi_i (A)}  =\sum_{\overline{A}\in\mathfrak{\overline{B}}_{n,k}} \left| \overline{A} \right|  \prod_{i=0}^n \left[ (n-i)! \right]^{\psi_i (\overline{A})}$$
matrices from  $\Pi_n$, which have at least $k$ elements that are component-wise equal to the respective elements of $\pi$.

And since $(n-n)!=0!=1$  and  $[n-(n-1)]!=1!=1$, then we finally obtain formula (\ref{fl3gfgfg}).
\hfill $\Box$

\begin{lemma}\label{lmm2}
For every integer $n\ge 2$
$$q(n,0)=q(n,1)=(n!)^{2n} =\left| \Pi_n \right| =\left| \Sigma _{n^2} \right| .$$
\end{lemma}

Proof. Let $k=0$. Then ${\mathfrak{B}}_{n,0}$ contains only the matrix, all elements of which are equal to 0. So $|{\mathfrak{B}}_{n,0}|=1$ and if $ A \in {\mathfrak{B}}_{n,0}$ then $|\overline{A} |=1$, $\psi_0 (A) =2n$ and $\psi_i (A) =0$ when $i\ge 1$. Therefore $\displaystyle q(n,0)= \sum_{\overline{A}\in \overline{\mathfrak{B}}_{n,0} } |\overline{A} |  \prod_{i=0}^{n-2} \left[ \left( n-i\right) ! \right]^{\psi_i (\overline{A})} = 1\cdot \left[(n-0)!\right]^{2n} \prod_{i=1}^{n-2} \left[ \left( n-i\right) ! \right]^0 = (n!)^{2n}$.

When $k=1$, there are $n^2$ matrices $A\in \mathfrak{B}_{n,1}$. It is easy to see that $|\mathfrak{\overline{B}} |=n$ and for every $\overline{A}\in\overline{\mathfrak{B}}_{n,1}$, $|\overline{A}|=n$, $\psi_0 (\overline{A})=2(n-1)$, $\psi_1 (\overline{A})=2$ and $\psi_i (\overline{A})=0$ for $i>1$. Therefore $q(n,1)=n^2 [(n-0)!]^{2n-2} [(n-1)!]^2 =(n!)^{2n-2} (n!)^2 =(n!)^{2n}$.
\hfill $\Box$

\begin{thm}\label{gl6_th1-bg}
Let  $A\in \Sigma_{n^2}$. Then the number $\xi_n$ of all matrices $B\in \Sigma_{n^2}$ which are disjoint with $A$ does not depend on $A$ and is equal to
\begin{equation}\label{gl6_main}
\xi_{n} =  \sum_{\overline{A}\in \overline{\mathfrak{B}}_n ,\; \varepsilon (\overline{A})\ge 2}  \left( -1\right)^{\varepsilon (\overline{A} )} \left| \overline{A} \right| \prod_{i=0}^{n-2} \left[ \left( n-i\right) ! \right]^{\psi_i (\overline{A})}
\end{equation}
\end{thm}

Proof. Let $n\ge 2$  be an integer. Then applying Theorem \ref{l2fhgg}, Lemma \ref{l3dskmat}, Lemma \ref{lmm2} and the principle of inclusion and exclusion we obtain that the number $\xi_n$  of all matrices $B\in \Sigma_{n^2}$  which are disjoint with  $A$ is equal to
$$
\begin{array}{ccc}
\xi_n & =  & \displaystyle |\Pi_n | +\sum_{k=1}^{n^2} (-1)^k q(n,k)  \\
& = & \displaystyle (n!)^{2n} -(n!)^{2n}+\sum_{k=2}^{n^2} (-1)^k q(n,k)  \\
& = & \displaystyle \sum_{k=2}^{n^2} (-1)^k q(n,k),
\end{array}
$$
where the function $q(n,k)$  is calculated with the help of formula (\ref{fl3gfgfg}). Thus we obtain the proof to formula (\ref{gl6_main}).
\hfill $\Box$

\begin{corollary}\label{th2-gl6}
The cardinality $\eta_{n}$ of the set of all  disjoint non-ordered  pairs of $n^2 \times n^2$ S-permutation matrices is equal to
\begin{equation}\label{nonordereddisjointpair_gl6}
\eta_{n} =\frac{(n!)^{2n}}{2} \xi_n
\end{equation}
where $\xi_n$ is described using formula \ref{gl6_main}.
\end{corollary}

Proof. It follows directly from formula (\ref{fcrl2}) and having in mind that the ''disjoint'' relation is symmetric and antireflexive.
\hfill $\Box$

\begin{corollary}\label{th3_gl6}
The probability $p_n$ of two randomly generated $n^2 \times n^2$ S-permutation matrices to be disjoint is equal to
\begin{equation}\label{probbility_gl6}
p_n = \frac{\displaystyle \xi_n}{\displaystyle  \left( n! \right)^{2n} -1}  ,
\end{equation}
where $\xi_n$ is described using formula (\ref{gl6_main}).
\end{corollary}

Proof. Applying Corollary \ref{th2-gl6} and formula (\ref{fcrl2}), we obtain:

$$p_n= \frac{\displaystyle \eta_{n}}{\displaystyle {\left| \Sigma_{n^2} \right| \choose 2}} = \frac{\displaystyle \frac{(n!)^{2n}}{2} \xi_n}{\displaystyle \frac{\left( n! \right)^{2n} \left( \left( n! \right)^{2n} -1\right) }{2}} = \frac{\displaystyle \xi_n}{\displaystyle  \left( n! \right)^{2n} -1} .$$
\hfill $\Box$

\section{An algorithm for counting}
There is one to one correspondence between the representation of the integers in decimal and in binary notations. So a square binary $n\times n$ matrix can be represented using ordered $n$-tuple of nonnegative integers, which belong to the closed interval $[0,\; 2^n -1]$.  Let the integer $a\in [0,\; 2^n -1]$. Then $a$ is represented uniquely in the form:
$$a=\sum_{u=0}^{n-1} b_u (a) 2^u ,$$
where $b_u (a)\in \mathfrak{B} =\{ 0,1\}$, $u=0,1,\ldots ,n-1$. We assume that we have implemented an algorithm for calculating the functions $b_u (a)$ for every $u=0,1,\ldots ,n-1$ and for every $a\in [0,\; 2^n -1]$. For example, in the programming languages C ++ and Java, $b_u (a)$ can be calculated using the expression

\begin{center}
\verb"bu = (a & (1<<u))==0 ? 0 : 1"
\end{center}

Let $A\in {\mathfrak B}_n $. With $\rho (A)$ we will denote the ordered $n$-tuple
$$\rho (A)=\langle x_{1} ,x_{2} ,\ldots ,x_{n} \rangle ,$$
where $0\le x_{i} \le 2^n -1$, $i=1,2,\ldots n$ and $x_{i} $ is the integer written in binary notation with the help of the $i$-th row of $A$.

We consider the set:
$$\begin{array}{lll} {{\mathfrak R}_n } & {=} & {\left\{\langle x_{1} ,x_{2} ,\ldots ,x_{n} \rangle \; |\; 0\le x_{i} \le 2^n -1,\; i=1,2,\ldots n\right\}} \\ {} & {=} & {\left\{ \rho(A)\, |\; A\in {\mathfrak B}_n \right\}} \end{array}$$

Thus we define the mapping
$\rho: {\mathfrak B}_n \to {\mathfrak R}_n ,$
which is bijective and therefore
${\mathfrak B}_n \cong {\mathfrak R}_n .$

If $A\in \mathfrak{B}_n$ and $\rho (A)=\alpha \in \mathfrak{R}_n$, then by analogy we define the numerical characteristics of the element $\alpha\in\mathfrak{R}_n$: $r_k (\alpha )= r_k (A)$, $c_k (\alpha ) =c_k (A)$, $\psi_k (\alpha ) = r_k (\alpha )+c_k (\alpha) =\psi_k (A)$, $k=0,1,2,\ldots ,n$ and $\varepsilon (\alpha )=\varepsilon (A)$. We assume $|\alpha |=|\overline{A}|$, where $\overline{A} =\{ B\in \mathfrak{B}_n \; |\; B\sim A \}$.

Let $\alpha=\langle x_1 , x_2 ,\ldots , x_n \rangle \in \mathfrak{R}_n$ and let $s$ be the number of different elements in $\alpha =\langle x_1 , x_2 ,\ldots , x_n \rangle$. Then the set $X=\{ x_1 , x_2 ,\ldots , x_n \}$ can be divide into parts $$X=X_1 \cup X_2 \cup \cdots \cup X_s$$ such that for every $k\in [s]$ and every $i,j\in [n]$, $i\ne j$ the condition $x_i ,x_j \in X_k$ is satisfied if and only if $x_i =x_j$. We assume $$z_i =\left| X_i \right| ,\quad i=1,2, \ldots s.$$

It is easily seen that
$$\displaystyle |\alpha |=\frac{n!}{\displaystyle \prod_{i=1}^s z_i !} .$$

Let $$\mathfrak{\overline{R}}_n =\{ \langle x_1 , x_2 ,\ldots , x_n \rangle \; |\; 0\le x_1 \le x_2 \le \cdots \le x_n \le 2^n -1 \} \subset \mathfrak{R}_n . $$

It is easily seen that $\mathfrak{\overline{B}}_n \cong \mathfrak{\overline{R}}_n$, which gives the basis to construct the following algorithm for calculating $\xi_n$:

\begin{algorithm} \label{alg1}
Calculation of $\xi_n$.\\
begin\\
$\xi_n :=0$ ;

For every $\alpha =\langle x_1 x_2 ,\ldots , x_n \rangle \in \mathfrak{\overline{R}}_n$ do

\{





\hspace{0.5 cm} $s:=1$;

\hspace{0.5 cm} $\varepsilon (\alpha ):=0$;






\hspace{0.5 cm} For $i=1,2,\ldots ,n$ do

\hspace{0.5 cm} \{

\hspace{1 cm} $z_s := z_s +1$ ; 

\hspace{1 cm} $t=0$;

\hspace{1 cm} For $u=0,1,\ldots , n-1$ do

\hspace{1 cm} \{

\hspace{1.5 cm} $t:=t+b_u (x_i )$;

\hspace{1 cm} \}

\hspace{1 cm} $r_t (\alpha ):=r_t (\alpha )+1$;

\hspace{1 cm} $\varepsilon (\alpha ) := \varepsilon (\alpha )+t$;

\hspace{1 cm} If $i<n$ and $x_i <x_{i+1}$ then $s:=s+1$;

\hspace{0.5 cm} \}

\hspace{0.5 cm} If $\varepsilon (\alpha ) = 0$ or $\varepsilon (\alpha )=1$ then go to next $\alpha$;

\hspace{0.5 cm} For $u= 0,1,\ldots ,n-1$ do

\hspace{0.5 cm} \{

\hspace{1 cm} $t:=0$;

\hspace{1 cm} For $i=1,2,\ldots ,n$ do

\hspace{1 cm} \{

\hspace{1.5 cm} $t=t+b_u (x_i )$;

\hspace{1 cm} \}

\hspace{1 cm} $c_t (\alpha ) :=c_t (\alpha )+1$;

\hspace{0.5 cm} \}

\hspace{0.5 cm} For $k=0,1,\ldots , n$ do

\hspace{0.5 cm} \{

\hspace{1 cm} $\psi_k  (\alpha ):= r_k (\alpha ) +c_k (\alpha )$;

\hspace{0.5 cm} \}

\hspace{0.5 cm} $\displaystyle |\alpha |:=\frac{n!}{\displaystyle \prod_{i=1}^s z_i !}$;

\hspace{0.5 cm} $\displaystyle T(\alpha):= (-1)^{\varepsilon (\alpha )} |\alpha | \prod_{i=0}^{n-2} \left[ \left( n-i\right) ! \right]^{\psi_i (\alpha)}$;

$\xi_n :=\xi_n +T(\alpha)$;

\}\\
end.

\end{algorithm}

\section{Conclusion}
On the basis of algorithm \ref{alg1} with programming language Java, we made a computer program for calculating $\xi_n$, $\eta_n$ and $p_n$ and we received the following results:\\
$$\xi_2 = 7$$
$$\xi_3 = 17\; 972$$
$$\xi_4 = 41\; 685\; 061\; 617$$
$$\xi_5 =232\; 152\; 032\; 603\; 580\; 176\; 504$$
$$\xi_6 = 7\; 236\; 273\; 578\; 711\; 450\; 275\; 537\; 707\; 547\; 657\; 855$$

$$\eta_2 = 56$$
$$\eta_3 = 419\; 250\; 816$$
$$\eta_4 = 2\; 294\; 248\; 126\; 968\; 596\; 791\; 296$$
$$\eta_5=71\; 871\; 209\; 790\; 288\; 983\; 974\; 921\; 874\; 964\; 480\; 000\; 000\; 000$$
$$\eta_6 = 7\; 022\; 228\; 210\; 556\; 132\; 949\; 916\; 635\; 069\; 726\; 824\; 032\; 981\; 704\; 989\; 720\; 182\; 784 \; \cdot \; 10^{13} $$

$$p_2 = 0.4666666666666667$$
$$p_3 = 0.38521058836137606$$
$$p_4 = 0.3786958223051558$$
$$p_5 = 0.37493849344703684$$
$$p_6 = 0.3728421644517476$$

For n = 2 and n = 3, the results that we get here coincide with the calculations made by hand in \cite{Yordzhev20151}, where we used a graph theory approach.


\end{document}